\documentclass[11pt]{amsart}
\usepackage{float}
\usepackage{amsmath, amssymb}
\usepackage{upgreek}
\usepackage{booktabs} 
\usepackage{xcolor, url}
\usepackage[all]{xy}
\usepackage{graphicx}
\usepackage{bm}
\usepackage{enumerate}
\theoremstyle{plain}

\newcommand{\td}{\text{d}}
\theoremstyle{definition}

\numberwithin{equation}{section}
\newcommand{\tr}{\operatorname{Tr}}

\usepackage{subfig, tikz}
\usetikzlibrary{decorations.pathmorphing}
\usepackage{graphics}

\usepackage{hyperref} %added
\hypersetup{colorlinks=true, linkcolor=blue, citecolor = blue, urlcolor=blue, filecolor=magenta}

\newcommand{\be}{\begin{equation}}
\newcommand{\ee}{\end{equation} }

\usepackage{pdflscape} %for landscape pages

\definecolor{ktbgreen}{RGB}{0, 145, 0}

\newcommand{\ktbadd}[1]{{\color{ktbgreen} 1}}

\begin{document}
\title[Spectrum of the Laplacian on the Page metric]{Spectrum of the Laplacian on the Page metric}
\author{Robie A.~Hennigar}
\address{Centre for Particle Theory, Department of Mathematical Sciences, Durham University 
\\
Durham DH1 3LE, United Kingdom}
	\email{robie.a.hennigar@durham.ac.uk}
\author{Hari K. Kunduri}
\address{Department of Mathematics and Statistics and Department of Physics and Astronomy\\
		McMaster University\\
		Hamilton, ON Canada}
	\email{kundurih@mcmaster.ca}
\author{Kam To Billy Sievers}
\address{Department of Physics and Astronomy\\
		McMaster University\\
		Hamilton, ON Canada}
	\email{sieversktb@mcmaster.ca}
\author{Yiqing Wang}
\address{Department of Mathematics and Statistics\\
		McMaster University\\
		Hamilton, ON Canada}
	\email{wang46@mcmaster.ca}

%\date{\today}

\thanks{H.~K.~Kunduri acknowledges the support of the  NSERC Grant  RGPIN-2018-04887.}
\thanks{K.~T.~B.~Sievers acknowledges the support of a NSERC PGS-D Scholarship.}

\begin{abstract}
We numerically construct the spectrum of the Laplacian on Page's inhomogeneous Einstein metric on $\mathbb{CP}^2 \# \overline{\mathbb{CP}}^2$ by reducing the problem to a (singular) Sturm-Liouville problem in one dimension.  We perform a perturbative analysis based upon a closely related, exactly solvable problem that strongly supports our results. We also study the spectrum of the Lichnerowicz Laplacian on symmetric traceless transverse two-tensors.  The method relies on both the isometries of the Page metric and pseudospectral methods to numerically solve the resulting ODEs. 
\end{abstract}

\maketitle 
\section{Introduction} A natural problem in mathematical physics is to determine the spectrum of the Laplacian operator $-\Delta$ on a compact Riemannian manifold $(M,g)$. The problem is closely associated to constructing a natural basis of complete eigenfunctions for the Laplacian. In turn the eigenvalues and eigenfunctions reveal important geometric and physical properties of fields on $(M,g)$, such as the time decay rates of solutions to the heat equation or yielding the fundamental vibrational frequencies of waves propagating on $M$.  For the case of the unit two sphere $\mathbb{S}^2$  with its maximally symmetric round metric, the solution to the problem is well known and the eigenfunctions are simply the spherical harmonics. 

More generally, consider a compact, four-dimensional Einstein manifold $(M,g)$ satisfying $\text{Ric}(g) = \Lambda g$, $\Lambda >0$.  The classical Lichnerowicz-Obata lower bound \cite{Besse} on the first (non-zero) eigenvalue $\lambda_1$  is 
\begin{equation}
\lambda_1 \geq \frac{4 \Lambda}{3} ,
\end{equation} with equality holding if and only if $(M,g)$ is isometric to $\mathbb{S}^4$ equipped with its round Einstein metric. The proof relies on an application of the Bochner identity. In the present work, we consider the four-dimensional \emph{gravitational instanton} $\mathbb{CP}^2 \# \overline{ \mathbb{CP}}^2$ equipped with its inhomogeneous Page metric $g$ \cite{Page:1978vqj}. In recent years, there has been a renewed interest in gravitational instantons. These are geodesically complete (possibly non-compact) Riemannian manifolds satisfying the Einstein condition $\text{Ric}(g) = \Lambda g$.  Geometries of this type arise naturally within the context of Euclidean quantum gravity as positive-signature analogues of Lorentzian black holes and are closely tied to black hole thermodyanmics \cite{Gibbons:1994cg}.  Of particular interest are hyper-K\"ahler examples (these are necessarily Ricci-flat, that is $\Lambda =0$) and Kronheimer achieved a complete classification~\cite{Kron}.  More recently, existence and uniqueness results for toric, Ricci flat gravitational imstantons which approach $\mathbb{S}^1 \times \mathbb{R}^3$ with its flat metric have been obtained \cite{KL}.  In the compact case with $\Lambda >0$, Page constructed a so called `rotating gravitational instanton' which from the geometer's perspective, is an Einstein manifold $(M,g)$ where $M = \mathbb{CP}^2 \# \overline{\mathbb{CP}}^2$ (a one-point blow up of $\mathbb{CP}^2$) or equivalently, the non-trivial $\mathbb{S}^2$-bundle over $\mathbb{S}^2$.  The associated metric $g$ is inhomogeneous with isometry group $SU(2) \times U(1)$ acting on $M$ with three-dimensional principal orbits. 

Using the fact that the Page metric is conformally K\"ahler,  Hall and Murphy produced the upper bound \cite{HallMurphy}
\begin{equation}
\lambda_1 < 1.89 \Lambda.
\end{equation}   By a result of Cao, Hamilton, and Ilmanen \cite{Cao}, since $\lambda_1 < 2 \Lambda$,  this establishes that the Page metric is linearly unstable under Ricci flow and can be destabilized by conformal perturbations. This complements a previous analysis of Young, who numerically demonstrated the existence of a negative eigenvalue for the Lichnerowicz operator on symmetric transverse traceless symmetric two tensors \cite{Young}. Within the framework of Euclidean quantum gravity,  the negative eigenvalue typically indicates that the solution is thermodynamically unstable, that is, it is an unstable critical point of the Euclidean action functional (see, for example \cite{Reall}). However, negative modes can arise which are not connected with the standard conditions for local thermodynamic stability~\cite{Dias}.

The main results of this work are as follows. We employ the $SU(2) \times U(1)$ isometry of the Page metric to reduce the eigenvalue problem for the Laplacian to a singular Sturm-Liouville problem in one dimension. We then use a technique based upon  pseudospectral analysis to numerically obtain a stronger upper bound on $\lambda_1$.   We also consider the analogous spectral problem for symmetric rank-2 tensor fields, i.e. the Lichnerowicz equation.  We repeat our analysis and recover Young's result with improved precision, although our negative value differs slightly that that obtained in \cite{Young}.

\section{Spectrum of the Laplacian of the Page metric}
\subsection{The Page metric}
The local form of Page's metric is
\begin{equation}
g = S \left[\frac{\td x^2}{A(x)} + 4 \alpha^2 A(x) \left(\td \psi + \frac{\cos\theta}{2} \td \phi \right)^2 + B(x) (\td \theta^2 + \sin^2\theta \td\phi^2) \right]
\end{equation} where 
\begin{equation}
A(x) = \frac{(3 - \nu^2 - \nu^2 (1+ \nu^2)x^2) (1-x^2) } {1 - \nu^2 x^2}, \qquad B(x) =  \frac{1 - \nu^2 x^2}{3 + 6\nu^2 - \nu^4}, \end{equation} and \begin{equation} S = \frac{3(1 + \nu^2)}{\Lambda}, \qquad  \alpha = (2 (3 + \nu^2))^{-1}.
\end{equation} Here  $x \in (-1,1)$, $\theta \in (0,\pi)$ and $\nu$ is a real parameter. The metric extends to a globally smooth metric on the the non-trivial $\mathbb{S}^2$ bundle over $\mathbb{S}^2$ provided the following identifications are satisfied: 
\begin{equation}\label{identifications}
(\psi, \phi) \sim (\psi + \pi, \phi + 2\pi), \qquad (\psi, \phi) \sim (\psi + 2\pi, \phi).
\end{equation} In particular, the Killing field $\partial_\psi$ degenerates smoothly at  $x = \pm 1$, the poles of the $\mathbb{S}^2$ fibres. Surfaces of constant $x$ have $\mathbb{S}^3$ topology.  The above metric is Einstein, i.e. 
\begin{equation}
\text{Ric}(g) = \Lambda g
\end{equation} provided that $\nu$ satisfies the quartic equation
\begin{equation}
\nu^4 + 4\nu^3 - 6 \nu^2 + 12 \nu - 3 =0.
\end{equation} There is a unique positive real root given approximately by  $\nu \approx 0.281702$ which we will take in the following.  Note that $\det g = 4 S^4 \alpha^2 B(x) \sin\theta$.  We will consider the eigenvalue problem 
\begin{equation}\label{eqn:eigeneq1}
-\Delta u = \lambda u.
\end{equation} Before considering the general case, we will restrict to $SU(2) \times U(1)$-invariant eigenfunctions $u = u(x)$. In this case the eigenvalue problem reduces to 
\begin{equation}\label{SLODE}
\partial_x \left[ \sqrt{B(x)} A(x) \partial_x u(x)\right] = -  S \sqrt{B(x)}  \lambda u \, .
\end{equation}

\subsection{Separation of Variables}  We now show that the eigenvalue equation admits separable solutions and can in fact be reduced to a singular one-dimensional Sturm-Liovuille problem. The Page metric has underlying isometry group $SU(2) \times U(1)$. This is most easily seen be writing it in terms of right-invariant one-forms on $SU(2)$: 
\begin{equation}
g = S \left[\frac{\td x^2}{A(x)} +  \alpha^2 A(x) \sigma_3^2 + B(x) (\sigma_1^2 + \sigma_2^2) \right]
\end{equation} where 
\begin{equation}
\sigma_1 = \sin \tilde\psi \td \theta - \cos\tilde\psi \sin\theta \td \phi, \quad \sigma_2 = \cos\tilde\psi \td \theta + \sin \tilde\psi \sin\theta \td \phi, \quad \sigma_3 = \td \tilde \psi + \cos \theta \td \phi, 
\end{equation} and $\tilde \psi:= 2\psi$.  These right-invariant forms satisfy the Maurer-Cartan equations $\td \sigma_i = -\tfrac{1}{2} \epsilon_{ijk} \sigma_j \wedge \sigma_k$.  In local coordinates, the generators of the $SU(2)_R \times U(1)_L$ are 
\begin{equation}
\begin{aligned}
R_1 &= -\cot\theta \cos\phi \partial_\phi - \sin \phi \partial_\theta + \frac{\cos\phi}{\sin\theta} \partial_{\tilde \psi}, \\
R_2 & = -\cot\theta \sin \phi \partial_\phi + \cos\phi \partial_\theta + \frac{\sin\phi}{\sin\theta} \partial_{\tilde \psi}, \\
R_3& = \partial_\phi, \qquad L_3 = \partial_{\tilde \psi}.
\end{aligned}
\end{equation}
 We can exploit the symmetry to reduce the eigenvalue problem to a single equation of Sturm-Liouville type by using an appropriate separation of variables.  We express the Page metric as 
\begin{equation}
g = 4 S \left[\frac{\td x^2}{4 A(x)} +  \alpha^2 A(x) \left( \td \psi + \frac{\cos\theta}{2} \td \phi \right)^2 + B(x) \hat{g} \right]
\end{equation} where 
\begin{equation}
\hat{g} = \frac{1}{4} (\td \theta^2 + \sin^2\theta \td \phi^2)
\end{equation} which we recognize as the Fubini-Study metric on $\mathbb{S}^2 \simeq \mathbb{CP}^1$ normalized so that $\text{Ric}(\hat{g}) = 4 \hat{g}$. The one-form appearing in the metric 
\begin{equation}
\hat{A} = \frac{\cos\theta}{2} \td \phi
\end{equation} is locally defined on $\mathbb{S}^2$ and satisfies $\hat{J} = \td \hat{A}$ where $\hat{J}$ is the K\"ahler form on $(\mathbb{CP}^1, \hat{g})$.  Let
$$ D  := \nabla_{\mathbb{S}^2} - i n A $$ where $n \in \mathbb{Z}$. 
We have
\begin{equation}
\begin{aligned}
D^2 & = \hat{g}^{ij} D_{i} D_{j} = g^{i j}\left( (\nabla_{\mathbb{S}^2})_{i} - i n A_{i} \right) \left((\nabla_{\mathbb{S}^2})_{j} - i n A_{j} \right) \nonumber \\ &= \Delta_{\mathbb{S}^2} - 2 i n A_{i} \hat{g}^{ij} (\nabla_{\mathbb{S}^2})_{j} - in \mbox{div}_{\hat{g}} A - n^2\hat{g}^{ij}
A_{i} A_{j} 
\end{aligned}
\end{equation} Since $\mbox{div}_{\hat{g}} A = \hat{g}^{ij} (\nabla_{\mathbb{S}^2})_{i} A_{j} = 0 $, 
\begin{align}
D^2 &= \Delta_{\mathbb{S}^2} - 2 i n A_{i} \hat{g}^{ij} (\nabla_{\mathbb{S}^2})_{j} - n^2\hat{g}^{i j}
A_{i} A_{j} 
\end{align}
We now compute $D^2$ explicitly. The Laplacian on $(\mathbb{S}^2, \hat{g})$ is
\begin{align}
\Delta_{\mathbb{S}^2} & = \frac{4}{\sin \theta} \partial_\theta \left(\sin \theta \partial_\theta  \right) + \frac{4}{\sin^2 \theta} \partial^2_\phi \end{align}
and the remaining terms are
\begin{align}
- 2 i n A_{i} \hat{g}^{ij} (\nabla_{\mathbb{S}^2})_{j}  &= -4 in \frac{\cos \theta}{\sin^2 \theta} \partial_\phi \nonumber \\ 
n^2\hat{g}^{ij} A_{i} A_{j} &= n^2 \frac{4}{\sin^2 \theta} \frac{\cos^2 \theta}{4} = n^2 \cot^2 \theta 
\end{align} which gives
\begin{align}
D^2 = \frac{4}{\sin \theta} \partial_\theta (\sin \theta \partial_ \theta) + \frac{4}{\sin^2 \theta} \partial_\phi^2 - n^2 \cot^2 \theta - 4 in \frac{\cos \theta}{\sin^2 \theta} \partial_\phi
\end{align} The operator $D^2$ is the charged Laplacian on $\mathbb{S}^2$ and its spectrum has been analyzed in detail in the context of $U(1)$ monopoles~\cite{Wu}.  Its eigenfunctions $Y(\theta, \phi)$ (suppressing the eigenvalue labels) are similar to the standard spherical harmonics, and satisfy \cite{KLR} 
\begin{equation}\label{spectrum}
D^{2} Y(\theta,\phi) = -\mu Y(\theta,\phi) 
\end{equation} where
$\mu \geq 0$ are a discrete family of eigenvalues with corresponding eigenfunctions $Y(\theta,\phi)$. The values taken by $\mu$ are, 
\begin{align}
\mu = \ell(\ell+2) - n^2\,, \mbox{ where } \, \ell = 2 k + |n| \,\, \mbox{ with } k=0,1,2,3...
\label{mulrelation}
\end{align} We can express the operation of the Laplacian operator associated to the Page metric in terms of $D^2$.  We seek separable solutions to 
\begin{equation}
\nabla^2 \Phi = -\lambda \Phi
\end{equation} of the form 
\begin{equation}
\Phi(x,\psi, \theta,\phi) = u(x) Y(\theta,\phi) e^{i n \psi}
\end{equation} where the periodicity $\psi \sim \psi + 2\pi$ requires $n \in \mathbb{Z}$.  The inverse of the Page metric $g$ in the $(x,\psi,\theta,\phi)$ coordinate chart is
\begin{gather}
g^{xx} = \frac{A}{S}, \qquad g^{\psi\psi} = \frac{1}{4\alpha^2 SA} + \frac{\cot^2\theta}{4 S B}, \qquad g^{\phi\phi} = \frac{1}{SB \sin^2\theta} \\ g^{\psi\phi} = -\frac{\cos\theta}{ 2S B \sin^2\theta}, \quad g^{\theta\theta} = \frac{1}{SB}
\end{gather} Using this and the Laplacian 
\begin{equation}
\nabla^2 \Phi = \frac{1}{\sqrt{\det g}} \partial_a \left(\sqrt{\det g} g^{ab} \partial_b \Phi \right)
\end{equation} the eigenvalue equation becomes, after using the separation ansatz
\begin{equation}
\frac{1}{S\sqrt{B} u} \partial_x [ A \sqrt{B} \partial_x u] - \frac{n^2}{4 \alpha^2 SA} + \frac{D^2 Y}{ 4 S B Y} = -\lambda
\end{equation}  where $D$ is the charged Laplacian operator on $\mathbb{S}^2$ defined above. Therefore, assuming $Y = Y(\theta,\phi)$ is an eigenfunction with eigenvalue $\mu$ for the given $n$, we get
\begin{equation}
\frac{1}{S\sqrt{B(x)} u(x)} \partial_x [ A(x) \sqrt{B(x)} \partial_x u(x)] - \frac{n^2}{4 \alpha^2 SA} -\frac{\mu}{ 4 S B } = -\lambda.
\end{equation} This leads to the final Sturm-Liouville equation
\begin{equation}
-\frac{1}{S\sqrt{B(x)}} \partial_x [ A(x) \sqrt{B(x)} \partial_x u(x)] = \left[ \lambda  - \frac{n^2}{4 \alpha^2 SA} -\frac{\mu}{ 4 S B }\right] u(x) 
\end{equation} For given integers $n$ and $k = 0,1,2,3\ldots$ the problem is to compute the corresponding $\lambda_{n,k}$. When $n = k =0$ so $\Phi = u(x)$, then we recover \eqref{SLODE}.

\subsection{Numerical Results} We employ a pseudospectral method to obtain the spectrum and eigenfunctions of the Laplacian on the Page metric. The numerical techniques, which are used both for the Laplacian and the Lichnerowicz operator, are described in detail in Appendix~\ref{sec:numerical}.

\begin{table}[h]
\begin{center}
\begin{small}
\begin{tabular}{c c c c c c c} 
 \toprule
 \multicolumn{7  }{c}{Laplacian Modes for Page Metric} \\
 \toprule
 $n$ & $k$  &  \multicolumn{5  }{c}{$\lambda_{N, n, k}$} 
 \\
 \midrule
 0 & 0 & 0 & 1.85251690621979 & 5.55511072563301 & 11.1094400455677 & 18.5152564002121
 \\
 0 & 1 & 2.20151767004730 & 4.10276203556860 & 7.79344389858214 &  13.3447981226032 & 20.7496082458243
 \\
 0 & 2 & 6.60146891570063 & 8.60169298548521 & 12.2721538919512 & 
17.8162828594157 & 25.2187121937168
\\
0 & 3 & 13.1940099324264 & 15.3461590383427 & 18.9949969583202 &
24.5254618756278 & 31.9233748061870
\\
0 & 4 & 21.9711957439182 & 24.3313655667953 & 27.9667687200084 & 33.4747507281330 & 40.8648223524852
\\
0 & 5 & 32.9239378762816 & 35.5508266297053 & 39.1923215047910 & 44.6674541140311 & 52.0447236550953
\\
\midrule
1 & 0 & 2.47278822016932 & 6.18734634277926 & 11.7439462625547 &
19.1505577680840 & 28.4082181994025
\\
1 & 1 & 5.73906807970214 & 9.51467128106754 & 15.0836636236597 &
22.4944425882384 & 31.7540476120613
\\
1 & 2 & 11.1813564846012 & 15.0588205997487 & 20.6502040479781 & 28.0679899485775 & 37.3307470900938
\\
1 & 3 & 18.7975109393924 & 22.8177272764530 & 28.4439648251217 & 35.8717888662822 & 45.1387882112358
\\
1 & 4 & 28.5847672388837 & 32.7885432777040 & 38.4652765244333 & 45.9066209583400 & 55.1788225190168
\\
1 & 5 & 40.5399430657669 & 44.9676920308993 & 50.7142102040491 & 58.1734108867898 & 67.4516717029354
\\
\midrule
2 & 0 & 6.94652453044506 & 12.5202199248593 & 19.9335667849511 & 29.1945205531415 & 40.3053303721494
\\
2 & 1 & 11.2814688789241 & 16.9219591058403 & 24.3618420429838 & 33.6356939266227 & 44.7537334325406
\\
2 & 2 & 17.7829446732259 & 23.5234041148932 & 31.0040472507313 & 40.2975471945021 & 51.4265039715026
\\
2 & 3 & 26.4498790489774 & 32.3231827179056 & 39.8598871482778 & 49.1801723955234 & 60.3238333569773
\\
2 & 4 & 37.2809295090199 & 43.3195060187479 & 50.9288969504986 & 60.2836590539944 & 71.4459635587983
\\
2 & 5 & 50.2745492141218 & 56.5102037560829 & 64.2103920748752 & 73.6080655744492 & 84.7931753765638
\\
\midrule
3 & 0 & 13.4241844219621 & 20.8552310645830 & 30.1252550740114 &  41.2412803729093 & 54.2060280076711
\\
3 & 1 & 18.8290895787195 & 26.3302703190900 & 35.6361972981928 & 46.7727794285795 & 59.7503486646426
\\
3 & 2 & 26.3952910470792 & 33.9943427949498 & 43.3510953723220 & 54.5167665456664 & 67.5124244590979
\\
3 & 3 & 36.1221531265956 & 43.8464861709966 & 53.2695094597091 & 64.4731162598423 & 77.4922746065115
\\
3 & 4 & 48.0088991038854 & 55.8854916907227 & 65.3908495786521 &  76.6416477479800 & 89.6899128471207
\\
3 & 5 & 62.0546379038340 & 70.1099249987414 & 79.7143611505433 & 91.0221108556643 & 104.105339193625
\\
\midrule
4 & 0 & 21.9065497148548 & 31.1938975332111 & 42.3201617079715 &  55.2914162802011 & 70.1104824309433
\\
4 & 1 & 28.3820039161805 & 37.7416690740012 & 48.9105160916212 & 61.9087188119635 & 76.7458556250276
\\
4 & 2 & 37.0154340588078 & 46.4711865532896 & 57.6971488670932 & 70.7315480883281 & 85.5929353417001
\\
4 & 3 & 47.8064222452126 & 57.3817420264504 & 68.6796186645005 &  81.7596888010515 & 96.6516427545189
\\
4 & 4 & 60.7544677247899 & 70.4724701939208 & 81.8573657302850 & 94.9928610023333 & 109.921871446885
\\
4 & 5 & 75.8589996332480 & 85.7423609513503 & 97.2297081639786 & 110.430713461739 & 125.403482270419
\\
\midrule
5 & 0 & 32.3939134142300 & 43.5369081835810 &  56.5189918644816 & 71.3454284480494 & 88.0189753506404
\\
5 & 1 & 39.9402331081855 & 51.1570843755572 & 64.1867597239542 & 79.0454537449158 & 95.7418033076496
\\
5 & 2 & 49.6422352144079 & 60.9537154020938 & 74.0447821741518 &  88.9451646249274 & 105.670983944151
\\
5 & 3 & 61.4996250546362 & 72.9262601706846 & 86.0926550173941 &
101.044314396956 & 117.806385167223
\\
5 & 4 & 75.5120551935158 & 87.0740711446609 & 100.329882521542 &
115.342595763238 & 132.147840182808
\\
5 & 5 & 91.6791322445776 & 103.396402218180 & 116.755877098339 &
131.839637732572 & 148.695144422478
\\ \bottomrule
\end{tabular}
\end{small}
\end{center}
\caption{Eigenvalues of the scalar Laplacian on the Page metric obtained with $\mathcal{N} = 250$. The overtone number $N$ increases from left to right within each row, with the left-most entry corresponding to $N = 0$. Units are such that $\Lambda = 1$.}
\label{scalar_modes_evals}
\end{table}

Our numerics exhibits rapid, typically exponential convergence for the eigenvalues of the Laplacian. This allows us to quickly obtain eigenvalues accurate to very high precision. We display in Table~\ref{scalar_modes_evals} the eigenvalues for $0 \le n,k \le 5$ and in each case the fundamental mode along with the first four overtones. In particular, we note that the lowest mode $\lambda_{1,0,0} = 1.85251690621979$ satisfies the bound $\lambda_{1,0,0} < 1.89 \Lambda$ obtained by Hall and Murphy~\cite{HallMurphy}. 

\begin{figure}[t]
	\includegraphics[width=\textwidth]{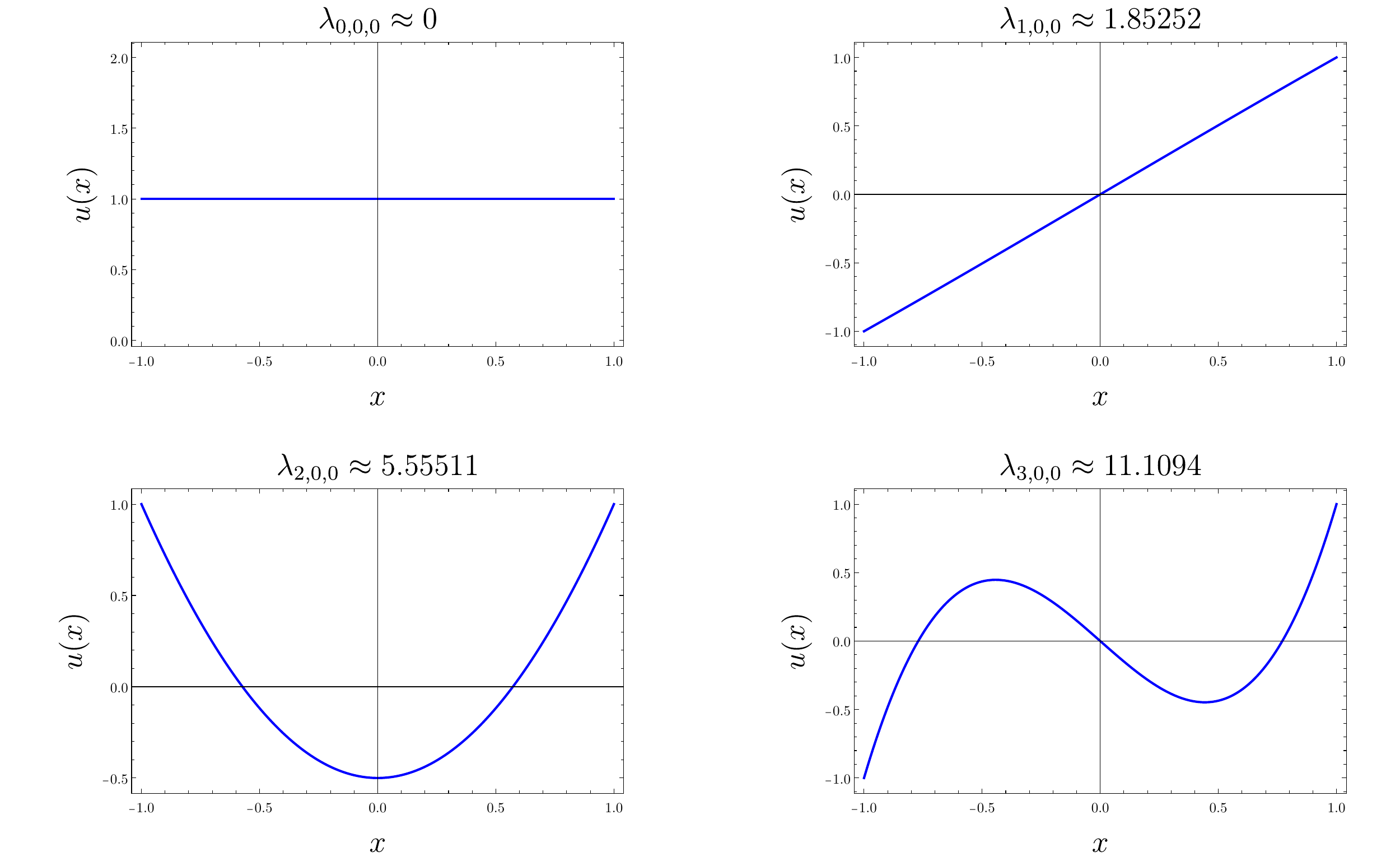}
	\caption{Plots of the first four eigenfunctions of the Laplacian with $n = k = 0$ scaled such that the right endpoint is at unity. Units are such that $\Lambda = 1$.}
	\label{fig:lap_eigenfuncs_00}
\end{figure}

Our method also allows us to obtain the eigenfunctions, of which we plot a selection in Figure~\ref{fig:lap_eigenfuncs_00}. The eigenfunctions bear a strong resemblance to the associated Legendre polynomials, though they are not exactly the same. We shall explain this similarly in the following section.  Additional plots of the eigenfunctions appear in Appendix~\ref{sec:EigenPlots}.

By examining the numerical results in the limit of large overtone number $N$, we are able to extract an asymptotic form that applies for all cases we have examined. We observe that 
\be 
\lambda_{N, n, k} \approx a N (N + 2n + 1) + \mathcal{O}(N^0) \, 
\ee
where we have determined $a \approx 0.925729$. We have arrived at this value of $a$ by numerically constructing the spectrum including over 500 overtones for several values of $(n,k)$ in the range $0 \le n,k \le 10$. We then fit that numerical data using {\sc Mathematica's} built-in ``Fit'' function, extracting the result given above. The coefficient $a$ is insensitive to the values of $(n,k)$, applying universally in all of the cases we have considered. We will provide analytical arguments in support of the generality of this result below.

\subsection{Approximation Methods}
We may obtain an approximate formula for the eigenvalues by considering the Page geometry as a perturbation of the exactly solvable $\nu =0$ eigenvalue problem. In this case, the geometry simplifies to 
\begin{equation}\label{homogeneous}
g = \frac{1}{\Lambda} \left[\frac{\td x^2}{1-x^2} + (1-x^2) \left( \td \psi + \frac{\cos\theta}{2} \td \phi \right)^2 + \td \theta^2 + \sin^2\theta \td \phi^2 \right].
\end{equation} As before, $x \in (-1,1)$ with $x= \pm 1$ corresponding to fixed points sets of $\partial_\psi$. The geometry extends globally to a smooth metric on the non-trivial $\mathbb{S}^2$ bundle over $\mathbb{S}^2$ provided the identifications \eqref{identifications} are imposed.  In contrast to the Page metric, this geometry is homogeneous with isometry group enhanced to $SO(3) \times SU(2)$. Here the first $SO(3)$ factor acts with two-dimensional orbits on the round sphere fibres parameterized by $(x,\psi)$. This additional symmetry allows for us to solve for the eigenvalues of the Laplacian in closed form.  The metric satisfies 
\begin{equation}
R^a_{~b} - \Lambda \delta^a_{~b} = \frac{(1-x^2)\Lambda}{8} \begin{pmatrix} 0 & 0 & 0 & 0 \\ 0 & 1 & 0 &\cos\theta \\ 0 & 0 & 1 & 0 \\ 0 & 0 & 0 & 1 \end{pmatrix}
\end{equation} where $x^a = \{x, \psi,\theta,\phi\}$. 

On the above metric, the equation for $u(x)$ takes a particularly transparent form. Directly, we have
\be 
(1-x^2) u''(x) - 2 x u'(x) + \left[\lambda^{(0)} - \frac{\mu}{4} - \frac{n^2}{1-x^2} \right] u(x) = 0 \, ,
\ee
where we introduced the notation $\lambda^{(0)} = \lim_{\nu \to 0} \lambda$. This is simply the associated Legendre differential equation, which allows one to read off the eigenfunctions (the associated Legendre polynomials) and the spectrum ($\lambda^{(0)} = \mu/4 + \ell(\ell+1)$) immediately. 

Given the simplicity of the problem when $\nu = 0$, our goal will be to obtain the spectrum for finite $\nu$ as a perturbation of the $\nu = 0$ result. We will begin by transforming the differential equation to Schr\"odinger form. Consider a Sturm-Liouville equation in standard form,
\be 
- \partial_x \left(p(x) \partial_x u(x) \right) + q(x) u(x) = \lambda w(x) u(x) \, .
\ee
For the  the Page metric, the functions can be readily identified,
\be 
p(x) = A(x) \sqrt{B(x)} \, , \quad q(x) = \left[\frac{n^2}{4 \alpha^2 S A(x)} + \frac{\mu}{4 S B(x)} \right] w(x) \, , \quad w(x) = S \sqrt{B(x)} \, .
\ee
We now introduce a new independent variable
\be 
s = \int_{-1}^x \sqrt{\frac{w(r)}{p(r)}} {\rm d}r 
\ee
and introduce a rescaled function
\be 
y(s) = \frac{u(s)}{m(s)} \, , \quad m(s) := \left(p w \right)^{-1/4} \, .
\ee
In these new variables, the original differential equation now takes the Schr\"odinger form 
\be 
\partial_s^2 y(s) = \left(V(s) - \lambda \right) y(s)
\ee
where the potential reads
\be 
V(s) = \frac{q(s)}{w(s)} + m(s)\frac{d^2}{ds^2} \left(\frac{1}{m(s)} \right) \, .
\ee

For the specific choice of $\nu$ required by the Page metric, the integral for the parameter $s$ cannot be carried out in closed form. However, for our purposes we only need this result perturbatively in $\nu$. We find:
\be 
x(s) = - \cos s + \frac{\nu^2  \sin s \left[\sin s\cos s - 3 s \right]}{6} + \mathcal{O}(\nu^4)\, .
\ee
This allows us to identify the potential and its leading correction $V = V_0 + \delta V + \mathcal{O}(\nu^4)$ where
\begin{align}
V_0 &= \frac{1}{4} \csc ^2(s) \left(\mu +4 n^2-(\mu -1) \cos ^2(s)-2\right) \, ,
\\
\delta V &= \frac{\nu^2 \csc^2(s)}{32} \bigg[5 \mu -4 \left(\mu +4 n^2+1\right) \cos (2 s)+32 n^2 s \cot (s)-16 n^2
\\
&-\mu  \cos (4 s)-8 s \cot (s)+12\bigg] \, . \nonumber
\end{align}

Let us first consider the solution to the unperturbed problem, that is with $\nu = 0$. In that case, we can easily obtain the spectrum and eigenfunctions exactly
\be 
y_{\ell, n}(s) = \sqrt{\frac{(1+2\ell) \left(\ell-n \right)!}{2 \left(n+\ell\right)!}} \sqrt{\sin s} P_\ell^n \left(\cos s \right) \, , \quad \lambda^{(0)} = \frac{\mu}{4} + \ell(\ell+1) \, ,
\ee
where $P_{\ell}^n(x)$ are the associated Legendre polynomials. Note that the quantum number $k$ does not actually enter into the explicit form of the eigenfunction, it appears in the eigenvalue only through its appearance in the parameter $\mu$. The eigenfunctions are normalized so that 
\be 
\int_{0}^\pi \left[y_{\ell, n}(s)\right]^2 {\rm d} s = 1 \, .
\ee

By the standard results of perturbation theory, we can obtain the correction to the spectrum by carrying out the following integral,
\be 
\lambda^{(1)} = \int_0^\pi \delta V(s) \left[y_{\ell, n}(s)\right]^2  {\rm d} s \, .
\ee
In evaluating this integral, we make use of the following intermediate results,
\begin{align}
\int_0^\pi \csc^2(s) \left[y_{\ell, n}(s)\right]^2  {\rm d} s &=  \frac{2\ell+1}{2 n} \, ,
\\
\int_0^\pi \cos(2s) \csc^2(s) \left[y_{\ell, n}(s)\right]^2  {\rm d} s &= \frac{2\ell + 1 - 4 n}{2 n} \, ,
\\
\int_0^\pi \cos(4s) \csc^2(s) \left[y_{\ell, n}(s)\right]^2 {\rm d} s &= \frac{32 n^3-32 n \ell ^2-32 n \ell +16 n+8 \ell ^3+12
   \ell ^2-2 \ell -3}{2 n (2 \ell -1) (2 \ell +3)} \,  ,
   \\
\int_0^\pi s \cot(s) \csc^2(s) \left[y_{\ell, n}(s)\right]^2  {\rm d} s &= \frac{(2 \ell +1) \left(4 n^2-4 n \ell -2 n-1\right)}{2
   n (2 n-1) (2 n+1)}\, .  
\end{align}
We then arrive at the following final result
\begin{align}
\lambda^{(1)} =\, \frac{\nu^2}{(2 n+2 N-1) (2 n+2 N+3)} \bigg[&4 N^4 +(16 n+8)
   N^3 
   \\
   &+    \left(-6 k^2-6 k (n+1)+20 n^2+21 n+1\right) N^2 
   \nonumber
   \\
   &+(2 n+1) \left(-6 k^2-6 k (n+1)+4 n^2+5 n-3\right) N
   \nonumber
   \\
   &-(2 n-1) \left(2
   k^2 (n+2)+2 k \left(n^2+3 n+2\right)-n (n+1)\right)\bigg] \nonumber \, ,
\end{align}
where we have replaced $\ell$ by the overtone number $N$, related by $\ell = n + N$.

One then expects that the spectrum of the Laplacian on the Page metric is  given by
\be 
\lambda = \lambda^{(0)} + \lambda^{(1)} + \mathcal{O}(\nu^4) \, .
\ee
upon substituting the particular value of $\nu \approx 0.281702$. The result of the perturbative calculation is in very good agreement with the numerical results, with the difference between the approximation and the ``exact'' numerical result arising in the second or third place after the decimal. This suggests that the perturbative expansion is not asymptotic but actually converges. We can use the approximation to obtain a simple form for the lowest lying modes ($N = 0$) of the scalar Laplacian on the Page metric. This reads, 
\be 
\lambda_{0, n, k} = \frac{3 n}{2} + n^2 + k n + k + k^2 + \left[\frac{\left(2 k^2 (n+2)+2 k \left(n^2+3 n+2\right)-n (n+1)\right)}{2
   n+3} \right] \nu^2 + \mathcal{O}(\nu^4) \, .
\ee 
We also note that while the intermediate integrals presented above require that $n > 0$, in the final result one can take $n = 0$ as well. Using the approximation to obtain the first nonzero eigenvalue for $n = 0$ and $k = 0$, we obtain
\be 
\lambda_{1, 0,0} = 2(1-\nu^2) + \mathcal{O}(\nu^4) \approx 1.84129 \, .
\ee
The above should be compared with our numerical result $\lambda_{1, 0,0}^{\rm num} \approx 1.85251690621979$. We see that the perturbative answer gives rather good agreement, with a relative error of approximately $0.6\%$. We have performed a broader comparison of the perturbative result to our numerics, considering $0 \le n,k \le 10$ and the first $250$ overtones for each choice $(n,k)$. Within this scanned parameter space, we find that the \textit{largest} relative error is $0.6\%$ for the $\lambda_{1, 0,0}$ mode.

On the other hand, we can expand our perturbative result in the limit of large overtone number,
\be 
\lambda_{N, n, k} = \left(1-\nu^2 \right)\left[N(1 + 2n + N) \right] + \mathcal{O}(N^0) \, .
\ee
The prefactor is approximately $1-\nu^2 \approx 0.920644$. Comparing this with our numerical results we find very good agreement. In fact, we find that the same functional form continues to hold
\be 
\lambda_{N, n, k}^{\rm num} = a \left[N(1 + 2n + N) \right] + \mathcal{O}(N^0) 
\ee
but with $a \approx 0.925729$. We expect an even more accurate approximation to this value of $a$ could be obtained by continuing the perturbative expansion to higher orders.

Continuing the perturbative analysis in \textit{general} is very tedious. However, if we are interested in a \textit{particular} eigenvalue, then it is not so difficult to obtain subleading corrections. For example, focusing on the $(n,k) = (0,0)$ mode, we obtain the following subleading corrections,
\be 
\lambda_{1,0,0} = 2 - 2 \nu^2 + \frac{202}{105} \nu^4 - \frac{362}{189} \nu^6 + \frac{24500302}{12733875} \nu^8 - \frac{105676778}{55180125} \nu^{10} + \mathcal{O}(\nu^{12})  \approx 1.852516461 \, .
\ee
This result agrees with the numerical result to six points after the decimal, with the agreement improving upon addition of further terms.

\section{Tensor modes} A much more involved problem is to determine the spectrum of the operator
\begin{equation}
\Delta_L(h)_{ab} = -(\nabla_b \nabla^b h_{ab} + 2 R_{acbd}h^{cd})~,
\end{equation} that is, solutions to
\begin{equation}\label{eqn:eigeneq2}
\Delta_L(h)_{ab} = \lambda h_{ab}~.
\end{equation} Here $h_{ab}$ is transverse, trace-free symmetric 2 tensor, which means it satisfies
\begin{equation}
 \nabla^a h_{ab} =0, \qquad \tr h:= g^{ab}h_{ab} = 0~.
\end{equation} These are gauge-fixing constraints that eliminate unphysical / redundant degrees of freedom. Even with these conditions, $h_{ab}$ naively $10 - 5 =  5$ independent components  and the eigenvalue problem will in general couple these together into a complicated 2nd order system. For the lowest eigenvalue, Young showed that one can reduce the system to a single master equation for $h_{00} = h_{00}(x)$:
\begin{equation}
h_{00} '' + C(x) h_{00}' + D(x) h_{00} = -\frac{\tilde \lambda} {A(x)} h_{00}
\end{equation} where
\begin{equation}\label{eqn:tildeLambda}
\lambda = \frac{\Lambda}{3 (1 + \nu^2)} \tilde \lambda
\end{equation} and
\begin{equation}\begin{aligned}
C(x) &= 3 E(x) + \frac{B'(x)}{B(x)} + \omega(x), \\ D(x) &= \frac{E(x)^2}{2} + \frac{3 E(x) B(x)'}{2 B(x)} - \left(\frac{B'(x)}{B(x)}\right)^2 + \frac{2b + 12 c x^2}{W(x)} + \frac{3 E(x) \omega(x)}{2} .
\end{aligned}
\end{equation} and
\begin{gather}
d = \frac{3 - \nu^2}{3 + 6 \nu^2 - \nu^4}, \quad b = -\frac{3 + \nu^4}{3 + 6\nu^2 - \nu^4}, \quad c= \frac{\nu^2(1 + \nu^2)}{3 + 6\nu^2 - \nu^4} \\
W(x)= c x^4 + b x^2 + d , \quad E(x) = \frac{2x (b + 2 c x^2)}{W(x)} \\
F(x) = \frac{2 x (3 + \nu^2 + x^2)}{B(x) (W(x))} \left[\frac{-1 + 2\nu^2 - \nu^4}{(3 + 6\nu^2 - \nu^4)^2}\right] \\
\omega(x) = \frac{1}{F(x)}\left[ \frac{4 B''(x)}{B(x)} + E(x)^2 - 2 \left(\frac{2 b + 12 c x^2}{W(x)}\right)\right]
\end{gather} 
Unfortunately, we have not found a value of $\nu$ for which the above equation admits exact solutions (it remains of Heun type). Therefore, we are unable to apply the perturbative techniques that were useful in obtaining analytical approximations to the spectrum of the Laplacian. We shall proceed via numerical techniques alone.

As explained by Young, an eigenfunction $h_{00}$ is normalizable (i.e. has finite $L^2$ energy) provided its singularity structure at $x = \pm 1$ is \textit{no stronger} than
\begin{equation} 
h_{00} \sim \frac{1}{1+x} \,  \quad \text{as} \quad x \to -1 \quad \text{and} \quad h_{00} \sim \frac{1}{1-x} \quad \text{as} \quad x \to 1 \, .
\end{equation}
Let us then examine the singularity structure of the differential equation. We begin with a Frobenius analysis near the endpoints of the domain. First, expanding near $x = -1$ we take the Frobenius ansatz
\be 
h_{00}(x) = \left(1-x^2\right)^{-p} \, \sum_{i=0}^n a_i (1+x)^i \, .
\ee
Substituting this expression into the differential equation and isolating the leading contribution allows us to determine two possible values for $p$:
\be 
p_\pm = \frac{11 + 3 \nu^2 \pm \sqrt{17 + 8 \nu^2 + \nu^4}}{2 \left(4 + \nu^2 \right)} \, .
\ee
Evaluating the above, we see that $p_- \approx 0.862628$ while $p_+ \approx 1.89224$. Since we demand that the solution be normalizable we must discard the ``$+$'' branch of the Frobenius solution. We repeat the same analysis at the other end point, $x = +1$, but it yields no further information.

Let us then define a new function $u(x)$ as
\be 
h_{00}(x) = \frac{u(x)}{\left(1-x^2\right)^{p_-}} \, .
\ee 
The leading singularity has been stripped off, and therefore the only requirement on $u(x)$ is that it should be regular at the endpoints of the domain. Solving the differential equation for $u(x)$ in a series expansion, the regularity requirement translates to two Robin boundary conditions:
\begin{align}
u'(\pm 1) =&\, \pm r(\nu, \lambda) u(\pm 1) \, ,
\\
r(\nu, \lambda) =& \,- \frac{ \left(-35732 \nu^3+53792 \nu^2-81660 \nu+19596\right) p_-  +50620 \nu^3-76217 \nu^2+115668 \nu-27783}{ (4 + \nu^2)(-3 + 2 \nu^2 + \nu^4) \left(15 + 4 \nu^2 - 2 p_- (4 + \nu^2) \right)}
\nonumber\\
&+ \frac{ \lambda \left(71 - 360 \nu + 233 \nu^2 - 152 \nu^3 \right)}{2(4 + \nu^2)(-3 + 2 \nu^2 + \nu^4) \left(15 + 4 \nu^2 - 2 p_- (4 + \nu^2) \right)} \, .
\end{align}

\begin{table}[t]
\begin{center}
\begin{small}
\begin{tabular}{c c c c c c} 
 \toprule
 \multicolumn{6  }{c}{Tensor Modes for Page Metric} \\
 \toprule
 $N$ &  $ \tilde\lambda_N$ & $N$ & $\tilde\lambda_N$ & $N$ & $\tilde\lambda_N$
 \\
 \midrule
 0 & $-$6.44142027579817 & 20 & 1379.23174426138 & 40 & 5156.63340079955
 \\
 1 & 12.7327722504243 & 21 & 1518.10938856071 & 41 & 5421.59387759233
 \\
 2 & 29.9361177572983 & 22 & 1649.05921793637 & 42  & 5666.26681454328
 \\
 3 & 55.3281097111635 & 23 & 1800.54515458952 & 43 & 5943.83557036251
 \\
 4 & 83.9221564235853 & 24 & 1942.86730153093 & 44 & 6199.88081482476
 \\
 \midrule
 5 & 121.930424598788 & 25 & 2106.96152640848 & 45 & 6490.05784947404
 \\
 6 & 161.900046603716 & 26 & 2260.65598788528 & 46 & 6757.47540128933
 \\
 7 & 212.518474497127 & 27 & 2437.35849822892 & 47 & 7060.26071461149
 \\
 8 & 263.861367308561 & 28 & 2602.42527218457 & 48 & 7339.05057365254
 \\
 9 & 327.088729504706 & 29 & 2791.73606609498 & 49 & 7654.44416552057
 \\
 \midrule
 10 & 389.804311324078 & 30 & 2968.17515109405 & 50 & 7944.60633168410
 \\
 11 & 465.640238602128 & 31 & 3170.09422722880 & 51 & 8272.60820199450
 \\
 12 & 539.728308281554 & 32  & 3357.90562224473 & 52 & 8574.14267519595
 \\
 13 & 628.172658686334 & 33 & 3572.43297963325 & 53 & 8914.75282386366
 \\
 14 & 713.633131796995 & 34 & 3771.61668391618 & 54 & 9227.65960403332
 \\
 \midrule
 15 & 814.685841339039 & 35 &  3998.75232184269 & 55 & 9580.87803098794
 \\
 16 & 911.518677424127 & 36 & 4209.30833483454 & 56 & 9905.15711806787
 \\
 17 & 1025.17971382173 & 37 & 4449.05225276177 & 57 & 10270.9838232507
 \\
 18 & 1133.38489152836 & 38 & 4670.98057404029 & 58 & 10606.6352171925
 \\
 19 & 1259.65423706624 & 39 & 4923.33277155864 & 59  & 10985.0702005543
 \\
 \bottomrule
\end{tabular}
\end{small}
\end{center}
\caption{Tensor modes for the Page metric obtained using the pseudospectral method with $\mathcal{N} = 551$. Units are such that $\Lambda = 1$.}
\label{tensor_modes_evals}
\end{table}

Using the techniques outlined in Appendix~\ref{sec:numerical}, we obtain the spectrum and eigenfunctions using a pseudospectral method. The first 60 eigenvalues are tabulated in Table~\ref{tensor_modes_evals}. Here we present the results with 15 digits of precision. It is possible to obtain significantly better precision than this (e.g. easily over 200 digits of precision for the low lying eigenvalues). In particular, we confirm the existence of a single negative mode, however obtaining a somewhat different answer than that quoted by Young, who reported $-5.75 < \tilde\lambda_0 < -5.74$ \cite{Young}.

\begin{figure}[H]
	\includegraphics[width=1\linewidth]{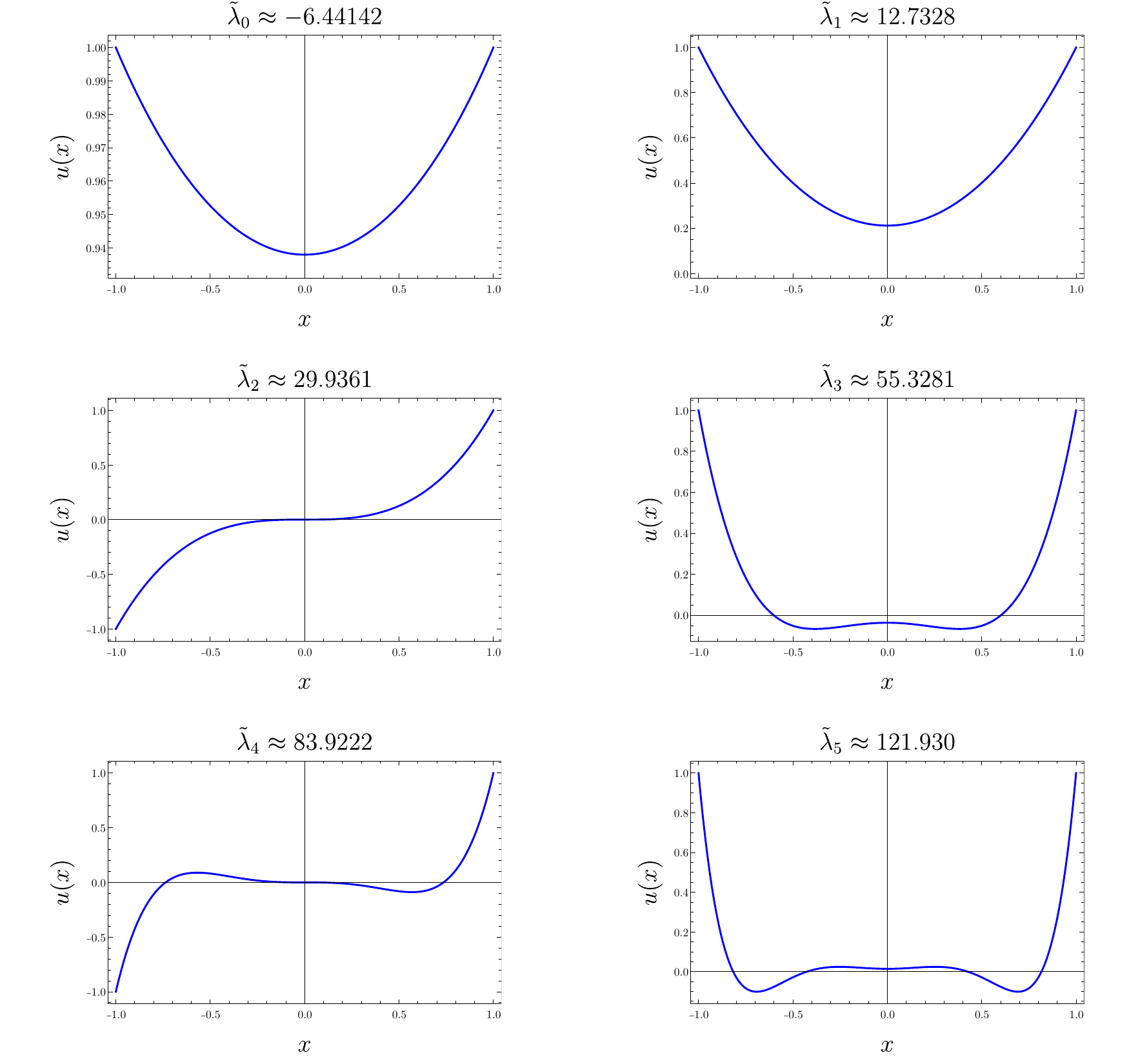}
	\caption{Plots of the $u(x)$ part of each associated eigenfunction in its domain $x\in(-1,1)$. The $u(x)$ functions and $\tilde\lambda$ values were taken from the $\mathcal{N}=300$ evaluation with the plotted functions scaled such that the right endpoint is at unity. Units are such that $\Lambda = 1$.}
	\label{fig:eigenfuncs}
\end{figure}

\appendix 

\raggedbottom
\section{Numerical Methods}
\label{sec:numerical}

\subsection{Pseudospectral Method}
Our primary means for obtaining the spectra and eigenfunctions of the differential operators is a pseudospectral method, which we review here. In both the scalar and tensor case, the problem reduces to a single ordinary differential equation on the interval $x \in [-1, 1]$. Schematically, the equations we consider all take the form
\be 
\tilde{\mathbf{H}}(h) = \lambda  \tilde{\mathbf{W}} h \, ,
\ee
where $\tilde{\mathbf{H}}$ is a second-order differential operator, $\tilde{\mathbf{W}}$ is a weight function, and $\lambda$ is the to-be-determined eigenvalue. We will describe the general features of our method. 

Given a differential equation for an undetermined function $h(x)$ on this interval, we first identify the singular behaviour of the solution near the endpoints using a Frobenius analysis. We select the branch of the Frobenius solution which has finite ${\rm L}_2$ energy and define a new function $u(x)$ according to the definition 
\be 
u(x) := (1+x)^{p_L} (1-x)^{p_R} h(x) \, ,
\ee
where $p_{L/R}$ is the degree of divergence of $h(x)$ near the left and right endpoints of the domain, respectively. Defined in this way, the only requirement on $u(x)$ is that it be regular. In terms of this new function, we can rewrite the differential equation as
\be 
\mathbf{H}(u) = \lambda \mathbf{W} u \, ,
\ee
where $\mathbf{H}$ and $\mathbf{W}$ differ from their earlier form by easily calculable factors.

We then solve the differential equation near each endpoint in a series expansion. At the left and right endpoints, the requirement of regularity results in a Robin boundary condition for the function $u(x)$ of the schematic character
\be 
u'(-1) = L^0(\nu) u(-1) + \lambda L^\lambda(\nu) u(-1) \quad \text{and} \quad u'(1) = R^0(\nu) u(1) + \lambda R^\lambda(\nu) u(1) \, ,
\ee
where $L^i(\nu)$ and $R^i(\nu)$ are constants that are determined from the perturbative solution. We have illustrated explicitly here that the boundary condition involves the eigenvalue directly in our case, but only linearly. 

To solve this problem we introduce a set of $\mathcal{N} + 1$ collocation points on the interval $[-1, 1]$. We take the Gauss-Lobatto points,
\be
x_i = \cos \frac{i \pi}{\mathcal{N}}  \quad \text{for} \quad i = 0, \dots, \mathcal{N} \, .
\ee
We then discretize the differential equation. The eigenfunction $u(x)$ becomes a vector that is defined at the Gauss-Lobatto points $u_i := u(x_i)$. For the derivatives appearing in the differential operator $\mathbf{H}$, we use the Chebyshev differentiation matrices, $D_{i,j}$. The explicit form of these objects can be found in, for example, section 2.4.2 of ref.~\cite{canuto2007spectral}. Similarly, we discretize the weight function by evaluating it at the Gauss-Lobatto points. The result is the following generalized eigenvalue problem
\be 
H_{i, j} u_j = \lambda W_{i,j} u_j \, 
\ee
where summation over $j$ is implied. To implement the Robin boundary conditions, the first row of the above matrix equation is replaced with the condition holding at $x = 1$, while the last row is replaced by the condition holding at $x = -1$. The terms that do not have factors of $\lambda$ are added to $H_{i,j}$, while those that do have factors of $\lambda$ are added to $W_{i,j}$. After the discretized equation has been obtained, we use the built-in ``Eigensystem'' function of {\sc Mathematica} to obtain the eigenvalues and eigenvectors of the system.  

\subsection{Convergence analysis}

\begin{figure}[H]
	\includegraphics[width=0.45\textwidth]{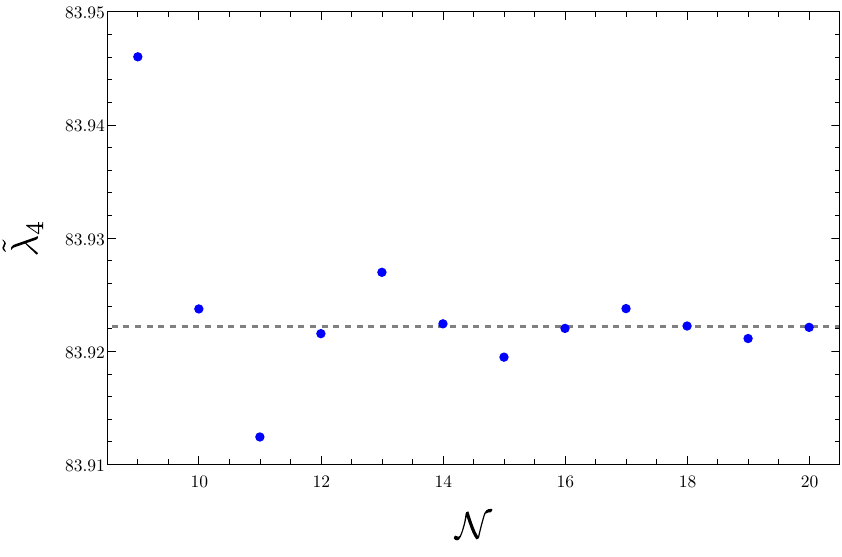}
	\quad
	\includegraphics[width=0.45\textwidth]{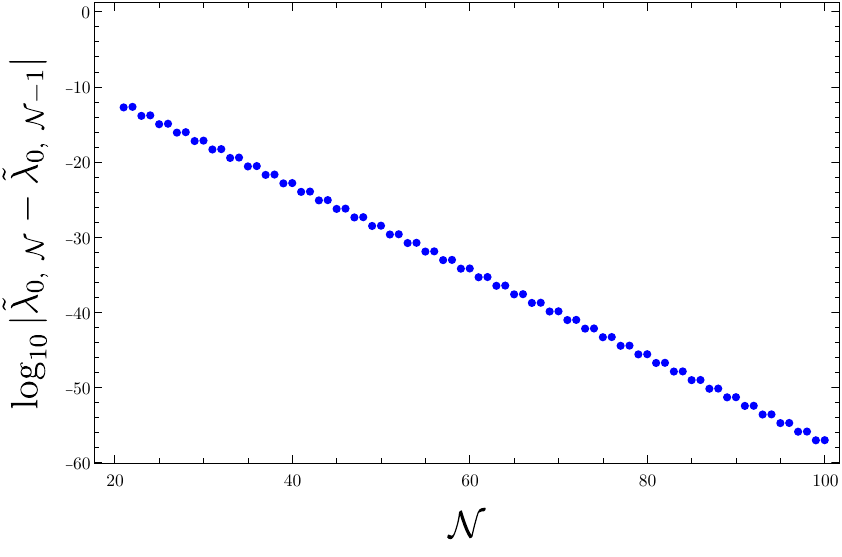}
	\caption{Left: Plot of the $N = 4$ eigenvalue ($\tilde{\lambda}_4$) of the Lichnerowicz operator on the Page metric for increasing $\mathcal{N}$. Right: Logarithmic plot of the difference between successive eigenvalue ($\tilde{\lambda}_0$) of the Lichnerowicz operator on the Page metric as a function of increasing $\mathcal{N}$. The choice of displaying the $\tilde\lambda_4$ and $\tilde\lambda_0$ cases were made as they best demonstrated the convergence behaviour --- other cases converged too quickly to showcase the convergence as seen on the left plot.}
	\label{fig:lich_converge}
\end{figure}

To assess whether a particular eigenvalue/vector combination is legitimate, we repeat the above discretization process for several values of $\mathcal{N}$ and keep only those pairs which exhibit convergence properties as a function of increasing $\mathcal{N}$. The eigenvectors obtained in this way provide a discrete approximation to the eigenfunctions of the differential operator, evaluated at the Gauss-Lobatto points. In all cases, the eigenvalues exhibit rapid convergence as $\mathcal{N}$ is increased, usually converging exponentially. We illustrate this for a particular eigenvalue in Fig.~\ref{fig:lich_converge}. 

In addition to the spectral method, we perform consistency checks using a more rudimentary shooting method. The technique is the same as that described in Appendix A of~\cite{Durgut:2022xzw} and therefore we do not describe it here. The results between the two methods are in agreement.

\raggedbottom
\pagebreak
\section{Further Plots of Laplacian Eigenfunctions}
\label{sec:EigenPlots}

\texttt{\begin{figure}[H]
			\includegraphics[width=\textwidth]{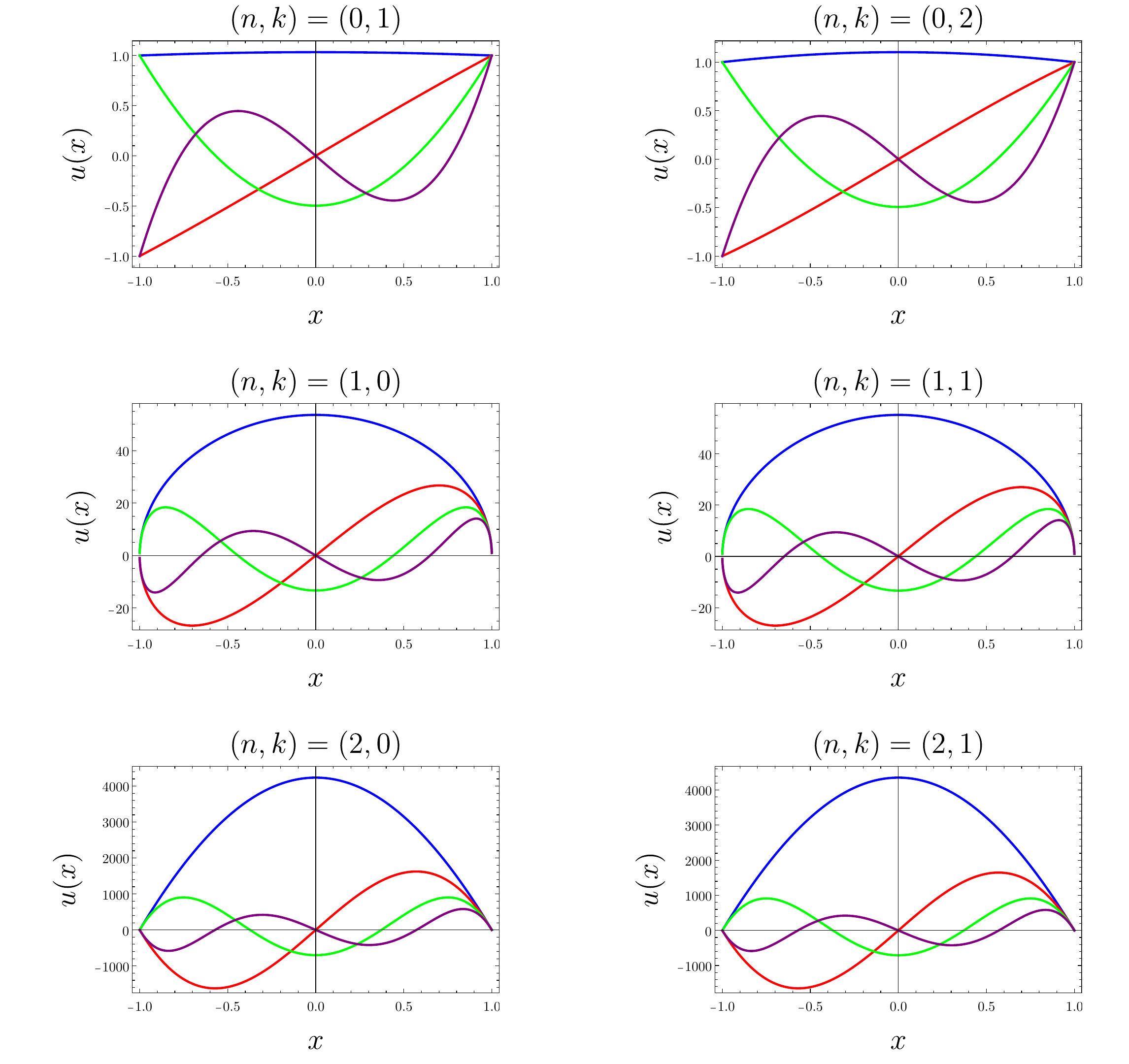}
			\caption{Plots of eigenfunctions of the Laplacian on the Page metric for six combinations of $n$ and $k$ with each color corresponding to varying overtone numbers (the blue, red, green, and purple colors represent the $N=0,1,2,3$ eigenfunctions). The eigenfunctions are scaled such that the right endpoint is at unity. Units are such that $\Lambda = 1$.}
			\label{fig:lap_eigenfuncs_01}
		\end{figure}}

In Figure~\ref{fig:lap_eigenfuncs_01} we present additional plots of the eigenfunctions for a selection of quantum numbers and overtones.		
		
%\texttt{\begin{figure}[t]
%	\includegraphics[width=\textwidth]{Billys_files/n0k1_4.pdf}
%	\caption{$\hdots$}
%	\label{fig:lap_eigenfuncs_01}
%\end{figure}
%\begin{figure}[h]
%	\includegraphics[width=\textwidth]{Billys_files/n1k0_4.pdf}
%	\caption{$\hdots$}
%	\label{fig:lap_eigenfuncs_10}
%\end{figure}
%\begin{figure}[h]
%	\includegraphics[width=\textwidth]{Billys_files/n1k1_4.pdf}
%	\caption{$\hdots$ .}
%	\label{fig:lap_eigenfuncs_11}
%\end{figure}

\end{document}